\documentclass[11pt, reqno,fleqn]{amsart}

\usepackage[margin=28mm]{geometry}

\usepackage[utf8]{inputenc}

\usepackage{amsmath}
\usepackage{amssymb,booktabs}
\usepackage{amsfonts}
\usepackage{graphicx}
\usepackage{amsthm,longtable}
\usepackage{enumerate, mathtools}
\usepackage{lscape}
\usepackage{dsfont, url}
\usepackage{xcolor, colortbl}
\usepackage{mathtools}

\usepackage{setspace}

\definecolor{Gray}{gray}{0.975}

\newtheorem{theor}{Theorem}[section]

\newtheorem{dfn}[theor]{Definition}

\newtheorem{prop}[theor]{Proposition}

\newtheorem{lem}[theor]{Lemma}
\newtheorem{guess}{Conjecture}
\newtheorem{problem}{Problem}
\newtheorem{Example}[theor]{Example}

\newtheorem{rem}[theor]{Remark}

\newcommand{\Z}{\mathds{Z}}
\newcommand{\Q}{\mathds{Q}}
\newcommand{\CP}{\mathds{C}\mathds{P}}

\newcommand{\T}{\mathds{T}}

\newcommand{\de}{\partial}
\newcommand{\dd}{\mathrm{D}}

\newcommand{\D}{\mathds{D}}

\newcommand{\inv}{$\mathds{T}^n$-invariant}

\renewcommand{\P}{\mathcal{P}}

\newcommand{\R}{\mathds{R}}

\newcommand{\C}{\mathds{C}}

\newcommand{\N}{\mathds{N}}

\newcommand{\K }{K\"{a}hler\ }
\newcommand{\KE }{K\"{a}hler-Einstein\ }
\newcommand{\I}{\mathcal{I}}

\begin{document}
\title[$\mathds{T}^n$-invariant K\"{a}hler-Einstein manifolds immersed in $\CP^n$]{$\mathds{T}^n$-invariant K\"{a}hler-Einstein manifolds immersed in \mbox{complex projective spaces}}

\author{Gianni Manno}
\address{(G. Manno) Dipartimento di Scienze Matematiche ``G. L. Lagrange'', Politecnico di Torino\\Corso Duca degli Abruzzi 24, 10129 Torino}
\email{giovanni.manno@polito.it}

\author{Filippo Salis}
\address{(F. Salis) Dipartimento di Scienze Matematiche ``G. L. Lagrange'', Politecnico di Torino\\Corso Duca degli Abruzzi 24, 10129 Torino }
\email{filippo.salis@polito.it}

\thanks{The first author gratefully acknowledges support by the projects ``Finanziamento alla Ricerca'' under the contract numbers 53\_RBA17MANGIO, 53\_RBA21MANGIO, and by the PRIN project 2022 “Real and Complex Manifolds: Geometry and Holomorphic Dynamics” (code 2022AP8HZ9). 
The second author is supported by the ``Starting Grant'' under the contact number 53\_RSG22SALFIL. Both authors are members of the GNSAGA of the INdAM}

\subjclass[2010]{53C55; 32Q20; 32Q40; 35J96; 35C11}
\keywords{\KE manifolds; \K immersions; \K toric manifolds; complex projective spaces; Calabi's diastasis function.}


\begin{abstract}
We give a complete list, for $n\leq 6$, of non-isometric \inv\ \KE manifolds immersed in a finite dimensional complex projective space endowed with the Fubini-Study metric. This solves, in the aforementioned case, a classical and long-staying problem addressed among others by Calabi and Chern.
\end{abstract}

\maketitle


\section{Introduction}

\subsection{Description of the problem}
Holomorphic and isometric immersions (from now on \emph{\K immersions}) into complex space forms, i.e. Kähler manifolds with constant holomorphic sectional curvature (see also Example \ref{cplxsf} below), have been  investigated starting from S. Bochner’s work \cite{bochner}. The general problem can be stated as follows:
\begin{problem}\label{genprob}
To determine whether a \K manifold can be \K immersed into a complex space form.
\end{problem}
Despite E. Calabi found in \cite{Cal} some criteria that allow, at least from a theoretical point of view, to treat the  Problem \ref{genprob}, we are in fact far from having an overall classification.
 Even in special cases of great interest, such as the \KE manifolds\footnote{Problem \ref{genprob} has also been  studied in more general context of K\"{a}hler-Ricci solitons, see  e.g. \cite{lmossa, pacific}.}, a complete solution of Problem \ref{genprob} is still lacking so far. Indeed, for \KE manifolds, Problem \ref{genprob} can be considered as solved only in the case of  immersions either into  hyperbolic or Euclidean spaces (see \cite{umehara2}), instead it  is still  open  for \KE manifolds  which are \K immersed into a complex projective space $\CP^N$. Partial results are contained,  e.g., in \cite{hano, hulinlambda, smyth} and also in  \cite{ch}, where S. S. Chern proved that the 1-codimensional \KE submanifolds of $\CP^N$ are either totally geodesics or the complex hyperquadrics, later extended by Tsukada \cite{ts} in the case of codimension 2.

The above considerations motivate the following definition.
\begin{dfn}
A \K manifold $(M,g)$ is called  \emph{projectively induced}, if $(M,g)$ can be \K immersed   into a (finite dimensional) complex projective space  $\CP^N$ endowed with the Fubini--Study metric $g_{FS}$,
namely the metric associated to the \K form  given in homogeneous coordinates by 
\begin{equation}\label{FS}
\frac{i}{2} \de\bar\de\log\left(|Z_0|^2+\mathellipsis+ |Z_N|^2\right).\end{equation}
\end{dfn}
Taking  into account that all the explicit examples of projectively induced \KE manifolds hitherto known are homogeneous manifolds (cfr. \cite{tak}), some conjectures have been put forward to justify this phenomenon (see e.g. \cite[Chap. 4]{loizedda}):
\begin{guess}\label{conj1}
A  $n$-dimensional complex manifold endowed with a projectively induced \KE  metric is an open subset of a complex flag manifold, i.e.  a compact simply-connected \K manifold on which its holomorphic isometry group acts   transitively.
\end{guess}
\begin{rem}\label{uniform}
If $n=1$, Conjecture \ref{conj1} is true. Indeed, by the Uniformization Theorem, the only $1$-dimensional \KE manifolds are open subsets of complex space forms. Moreover, the only complex space forms that admit a \K immersion into $(\CP^N,g_{FS})$ are complex projective spaces endowed with an integer multiple of the Fubini-Study metric (see \cite{Cal}).
\end{rem}

\subsection{Description of the main result}
In the present paper we focus our attention  to Conjecture \ref{conj1} in the case of \inv\ \KE manifolds  (which includes the remarkable class of the toric manifolds, see Example \ref{toricman} below).
In order to better describe our main result, we need some further definitions.

\begin{dfn}\label{inv}
A connected \K manifold $(M,g)$ is called  \inv\ if there exists an effective, biholomorphic and Hamiltonian action on $M$  of the $n$-dimensional real torus $\mathbb{T}^n$ with at least one fixed point.
\end{dfn}
The following examples provide some important classes of \inv\ \K manifolds.
\begin{Example}[Complex space forms] \label{cplxsf}
Let $z=(z_1,\mathellipsis, z_n)$ be complex holomorphic variables on $\C^n$. 
\begin{itemize}
\item Complex Euclidean spaces $$(\C^n,g_{euc}),$$ 
where $g_{euc}$ is the metric associated to the \K form $i\de\bar\de \|z\|^2$,
\item complex hyperbolic spaces
 $$\left(\D^n, g_{hyp}\right),$$ 
where $\D^n=\{z\in\C^n\ |\ \|z\|^2<1\}$ and $g_{hyp}$ is the metric associated to the \K form  $-\frac{i}{2}\de\bar\de\log(1-\|z\|^2)$ and
\item complex projective spaces 
$$\left(\CP^n, g_{FS}\right),$$ 
where $g_{FS}$   is the metric associated to the \K form $\frac{i}{2} \de\bar\de\log\left(1+\|z\|^2\right)$,
cfr. \eqref{FS},
\end{itemize}
are examples of \inv \ \K manifolds.
\end{Example}
\begin{Example}[Toric \K manifolds]\label{toricman}
Toric $n$-dimensional \K manifolds are  remarkable examples of \inv\ \K manifolds. We recall that the former  are defined as compact $n$-dimensional \K manifolds where it is defined an effective, biholomorphic and Hamiltonian action of the $n$-dimensional real torus $\T^n$.
As a consequence of the Atiyah and Guillemin-Sternberg's theorems (see \cite{atiyah, GS,guil}), a toric  $n$-dimensional  manifold has at least $n+1$ fixed points with respect to the $\T^n$-action. 
Complex projective spaces $\CP^n$ are basic examples of toric manifolds.
\end{Example}

As we said, the main object of the present paper will be \inv\  \KE manifolds. In particular, in this context,  Problem \ref{genprob} can be translated into the following one.
\begin{problem}\label{problem}
Classify all the projectively induced \inv\  \KE  manifolds.
\end{problem}
%
%
Since complex projective spaces are the only irreducible \inv\ flag manifolds and since only the integer multiples of the Fubini-Study metric  are projectively induced  (see \cite{Cal,loizedda}),  in the specific case of  \inv\ \K metrics,   Conjecture  \ref{conj1} reads as:
\begin{guess}\label{subconj}
The only projectively induced  \inv\ \KE manifolds are open subsets of $\CP^{n_1}\times\mathellipsis\times\CP^{n_k}$, with $n_1+\mathellipsis+n_k=n$, endowed with the \K metric
$$q\left(c_1 g_{FS} \oplus \mathellipsis \oplus c_k g_{FS}\right), $$
where $k$ and $q\in\Z^+$, $c_i=\frac{1}{G^{k-1}}\prod_{{ j\neq i}}(n_j+1)$ for $i=1,\mathellipsis,k$ and $G=\mathrm{gcd}(n_1+1,\mathellipsis, n_k+1)$, namely   the greatest common divisor between $n_1+1,\mathellipsis, n_k+1$.
\end{guess}
Our main goal consists in proving  Conjecture \ref{subconj} for $n\leq 6$, i.e., the following theorem.
\begin{theor}\label{main}
Conjecture \ref{subconj} is true for $n\leq 6$.
\end{theor}
We note that Conjecture \ref{subconj} has been already proved\footnote{We also note that Conjecture \ref{subconj} is true if we replace the assumption of $\mathds{T}^n$-invariance with the more restrictive one of $\mathds{U}(n)$-invariance  (cfr. the results achieved in  \cite{extr} thanks to some techniques based on the Calabi's diastasis function developed in \cite{mathz}).} in the case $n=2$ in \cite{ny} and in the case $n\leq 4$ for toric manifolds in \cite{loitoric}. In both papers the results have been achieved thanks to a careful analysis of a family of real Monge-Ampère equations, that required some demanding computations.
The approach we adopt in the present paper is slightly different: our main analytical object is a single Monge-Ampère equation whose we study the gradient map of its solutions. A central role will be played by the closure of the image of such maps, which in the context we are dealing with, turn out to be the so-called Delzant polytopes. Some further explanations are contained in Section \ref{strategy} below.


\subsection*{Notation and conventions}
If $I=(I_1,\mathellipsis,I_n)\in\N^n$ is a multi-index, then the length of $I$ is defined as $|I|:=\sum_{\alpha=1}^n I_\alpha$. If $w=(w_1, \mathellipsis,w_n)$, then $w^{I}$ denotes  the monomial in $n$ variables $\prod_{\alpha=1}^n w_\alpha^{I_\alpha}$.

\section{Proof of Theorem \ref{main}}

In Sections \ref{strategy}--\ref{simplex} we obtain results for arbitrary dimension $n$. In Sections \ref{1dim}--\ref{6dim} we used the results obtained in the previous sections in the case $n\leq 6$ to prove Theorem \ref{main}.
\subsection{Strategy of the proof and description of the paper}\label{strategy}
The first step for the proof of Theorem \ref{main} consists in constructing a bijective correspondence between \KE metrics object of  Problem \ref{problem} and polynomial solutions of a certain class of Monge-Ampère equations depending on the Einstein constant. More precisely,  we are going to describe in Section \ref{reduction} how such bijective correspondence is realized by means of a particular \K potential, namely the Calabi's diastasis function. In a suitable coordinate system $x=(x_1,\mathellipsis,x_n)$, after a suitable normalization of the Einstein constant (such normalization can be assumed without loosing of generality), the above bijective correspondence is obtained in terms of solutions of the
 $n$-dimensional Monge-Ampère equation 
\begin{equation}\label{mongeampere2}
 \det D^2u=e^{-u},\quad u(x)=u(x_1,\mathellipsis,x_n)
 \end{equation} of the type
\begin{equation}\label{type2}
u(x)=\log\sum_{I\in\I} a_I e^{I\cdot x}-\sum_{i=1}^n x_i,\end{equation}
where  $\I$ is a finite subset of $\N^n$ 
containing every multi-indices $I=(I_1,\mathellipsis,I_n)$ such that  $a_I\geq 0$ and  $a_{I}=1$ for $|I|\leq 1$. Furthermore, taking into account what we said above, Conjecture \ref{subconj} can be reformulated in simpler way, see Conjecture \ref{conj3}.

The second step is contained in Section \ref{gradient} and it concerns the study of the gradient map of the solutions to Monge-Ampère equation \eqref{mongeampere2}. In particular, by considering \cite{del} and \cite{mcduff} together with some technical lemmas, we  prove that the closure of the gradient image  has to be the dual of a smooth Fano polytope, in particular it is a reflexive polytope, having the origin as barycenter. Conversely, any convex polytope having the origin as barycenter is the closure of the  image of the gradient map of a convex solution to \eqref{mongeampere2} (see \cite{bb}), not necessary of type \eqref{type2}.  Section  \ref{gradient} ends by proving that decomposable convex polytopes are associated to solutions to the $n$-dimensional Monge-Ampère equation \eqref{mongeampere2} if and only if they are cartesian products of polytopes associated to solutions to lower dimensional Monge-Ampère equations of type \eqref{mongeampere2}.

The third step is contained in Section \ref{technical}. Taking into account the results of \cite{ny}, we find some constraints on the shape of polytopes associated to solutions of type \eqref{type2} to the Monge-Ampère equation \eqref{mongeampere2}.

Finally, in Section \ref{final},  we conclude the proof of Theorem \ref{main} through a case-by-case analysis  of the few polytopes satisfying all the properties  identified in the previous sections. In particular, we exploit the classification of reflexive polytopes by M. Øbro \cite{obro2} in order to prove that the only polytopes associated to solutions of the equation \eqref{mongeampere2} of type \eqref{type2} are either the simplex or a cartesian product of lower dimensional simplices.

\subsection{Reduction of Problem \ref{problem} to the existence of polynomial solutions of a $n$-dimensional Monge-Ampère equation}\label{diastasis}\label{reduction}

 Since \K potentials of the Fubini-Study metric  $g_{FS}$ are   real analytic functions, a \K potential $\Phi$ of a \K metric defined as holomorphic pullback of $g_{FS}$ is forced itself to be real analytic. Thus, considering a holomorphic system of coordinates on $U\subseteq \C^n$
$$z=(z_1,\mathellipsis,z_n),$$
 $\Phi$ is equal to  its  power expansion around the origin:
\begin{equation}\label{powexdiastc}
\Phi(z)=\sum_{I,J\in\N^n}a_{IJ}z^{I}\bar z^{J}.
\end{equation}
The function
\eqref{powexdiastc} can be complex analytically extended to a function $\tilde\Phi$ on a neighborhood of the diagonal in $U \times\overline{ U}$, where $\overline{ U}$ is the conjugate of $U$, thus defining
the \emph{diastasis function} $\dd_0:U\to\R$ for $g$:
$$\dd_0(z)=\tilde\Phi(z,\bar z)-\tilde\Phi(z,0)-\tilde\Phi(0,\bar z)+\tilde\Phi(0,0).$$

Moreover, for any K\"ahler manifold with real analytic metric, there exists a coordinates system, that we will still denote by $z=(z_1,\mathellipsis,z_n)$, in a neighbourhood of each point,  such that
\begin{equation}\label{bochnercoordinates}
\dd_0(z)=\sum_{\alpha=1}^n|z_\alpha|^2+\psi,
\end{equation}
where $\psi$ is a power series with degree $\geq 2$ in both $z$ and $\bar z$. 
\begin{dfn}
Coordinates $(z_1,\mathellipsis,z_n)$ for which \eqref{bochnercoordinates} holds true  are called \emph{Bochner's coordinates} for the metric $g$.
\end{dfn}
Bochner's coordinates are uniquely determined up to unitary transformations (cfr. \cite{bochner,Cal,hulin,hulinlambda,tian4}).
%

\begin{lem}\label{diastpolin}
Let $ (M,g)$ be a projectively induced  \inv\ \K manifold. Let  $z=(z_1,\mathellipsis,z_n)$  be Bochner's coordinates for $g$ centered at a fixed point $p$ of the toric action.
Then the diastasis function
$\dd_0(z)$  can be written  as
\begin{equation}\label{diastbochpi}
\dd_0(z)=\log\left(P(z)\right),
\end{equation}
where
\begin{equation}\label{eq.Pz}
P(z)=\sum_{I\in\I} a_I|z^{I}|^2
\end{equation}
with $a_I>0$ and  $a_I=1$ for every $|I|\leq 1$. Here $\I$ is a finite subset of $\N^n$.
\end{lem}
\begin{proof}
Let $Z_0,\dots, Z_N$ be homogeneous coordinates on $\CP^N$, with $N\geq n$. Let  $\zeta_j=\frac{Z_j}{Z_0}$ be the affine coordinates around the point $[1, 0\dots, 0]$ on $ U_0=\{Z_0\neq 0\}$. Let $f:M\to \CP^N$ be a \K immersion.
Up to a unitary transformation of $\CP^N$ and, if necessary, by restricting $f$ to an open neighborhood $V$ of $p$, we can assume that $f(p)=[1, 0\dots, 0]$ and $f(V)\subset U_0$.  

By \cite[Theorem 7]{Cal}, there exist Bochner's coordinates on $U_0$ such that the restriction of $f$ to an open neighborhood of $p\in M$ where are defined  Bochner's coordinates $z=(z_1,\mathellipsis,z_n)$ for $g$
can be written in coordinates as the graph of a holomorphic function:
$$ z=(z_1, \dots , z_n)\mapsto (z_1, \mathellipsis ,z_n, f_{n+1}(z), \dots , f_N(z)),$$
where
$$f_j(z)=\sum_{I\in\N^n}\alpha_{jI}z^{I},\quad \ j=n+1, \dots, N.$$
Since one can check that the affine coordinates $\zeta_j$ on $U_0$ are Bochner's coordinates for the Fubini--Study metric $g_{FS}$ and by considering also that the diastasis function is hereditary (see \cite[Prop. 6]{Cal}) in the sense that the diastasis function of the \K immersed manifold is given by the composition of the diastasis of the ambient space, that in our case  is $\log(1+\sum_{j=1}^N|\zeta_j|^2)$, with the \K immersion, we get
$$\dd_0(z)=\log\left(1+\sum_{j=1}^n|z_j|^2+\sum_{j=n+1}^N|f_j(z)|^2\right).$$
It is not hard to see that the Diastasis function centered at a fixed point of the toric action depends only on the moduli of the coordinates $z$ (see e.g. \cite{tinv}). Hence,
$f_j$'s needs to be monomials
in $z$ and formula \eqref{diastbochpi} follows.
\end{proof}


\begin{lem}\label{einsteinconst}
The Einstein constant of a projectively induced  \inv\ \KE manifold is a positive rational number.
\end{lem}
\begin{proof}
Let $z=(z_1,\mathellipsis,z_n)$ an arbitrary holomorphic coordinates system.
A K\"ahler metric $g$ with diastasis function $\dd_0(z)$  is  Einstein   if and only if  there exists $\lambda\in\R$ such that
$$\lambda\frac{i}{2} \de\bar\de \dd_0=-i\de\bar\de\log\det(g_{\alpha\bar\beta}).$$
Hence, by the $\partial\bar\partial$-lemma,  there exists a holomorphic function $\varphi$ such that
\begin{equation}\label{eq.complex.MAE}
\det(g_{\alpha\bar\beta})=e^{-\frac{\lambda}{2}(\dd_0+\varphi+\bar \varphi)}.
\end{equation}

Let now  set Bochner's coordinates, that we keep denoting by $z$.
 By comparing the series expansions of both sides of the equation \eqref{eq.complex.MAE}, we get that $\varphi+\bar \varphi$ is forced to be zero (cfr. e.g. \cite{hulinlambda,salis}).

In particular, since $\dd_0$ is the diastasis of a projectively induced \inv\ \K metric, by taking into account Lemma \ref{diastpolin}, the equation \eqref{eq.complex.MAE} reads as an equality between polynomials of real variables $$x_\alpha:=|z_\alpha|^2,$$ more precisely:
\begin{equation}\label{MA*}
\frac{\det\left[\left(P \frac{\partial^2 P}{\partial{x_\alpha}\partial{x_\beta}}-\frac{\partial P}{\partial{x_\alpha}}\frac{\partial P}{\partial{x_\beta}}\right)x_\alpha+P\frac{\partial P}{\partial{x_\alpha}} \delta_{\alpha\beta}\right]_{1\leq\alpha,\beta\leq n}}{P^{n-1}}=P^{-\frac{\lambda}{2}+n+1}.
\end{equation}
Hence, $\lambda$ needs to be a rational number. Moreover, by comparing the degrees of both sides of \eqref{MA*}, we get  $\lambda\geq 2\frac{n}{\deg P }>0$.\end{proof}

\begin{rem}\label{l1}
In view of Lemma \ref{einsteinconst}, we have
$\lambda=2 \frac{s}{q}\in\Q^+$, where we are assuming $\mathrm{gcd}(s,q)=1$. 
Let  $P(|z_1|^2,\mathellipsis,|z_n|^2)$ be a polynomial solution of type \eqref{eq.Pz}  to \eqref{MA*}. Since $\mathrm{gcd}(2nq,s)=1$, $P$ is forced to be the $q$-th power of a polynomial, namely
\begin{equation}\label{eq.Px}
P(|z_1|^2,\mathellipsis,|z_n|^2)=P(x_1,\mathellipsis,x_n)= \sum_{I\in\I} a_Ix^{I}, \quad \text{with } a_I=1 \text{  if  } \, |I|\leq 1,
\end{equation}
must be equal to
$$R\left(\frac{ x_1}{q},\mathellipsis,\frac{x_n}{q}\right)^q$$
where $R(x_1,\mathellipsis, x_n)$ is a polynomial of type \eqref{eq.Px}.
One can easily check that $R\left( |z_1|^2,\mathellipsis,|z_n|^2\right)=R\left( x_1,\mathellipsis,x_n\right)$ is a solution to the equation
\begin{equation}\label{MA**}
\frac{\det\left[\left(R \frac{\partial^2 R}{\partial{x_\alpha}\partial{x_\beta}}-\frac{\partial R}{\partial{x_\alpha}}\frac{\partial R}{\partial{x_\beta}}\right)x_\alpha+R\frac{\partial R}{\partial{x_\alpha}} \delta_{\alpha\beta}\right]_{1\leq\alpha,\beta\leq n}}{R^{n-1}}=R^{n+1-s}.
\end{equation}
Moreover, if a polynomial $R(x_1,\mathellipsis,x_n)$ of type \eqref{eq.Px} is a solution of \eqref{MA**} for some $s\in\N$, then 
$$P(x_1,\mathellipsis,x_n)=R\left(\frac{x_1}{s},\mathellipsis,\frac{x_n}{s}\right)^{s}$$
 is a polynomial solution of the same type to the equation 
\begin{equation}\label{MAeq}
\frac{\det\left[\left(P \frac{\partial^2 P}{\partial{x_\alpha}\partial{x_\beta}}-\frac{\partial P}{\partial{x_\alpha}}\frac{\partial P}{\partial{x_\beta}}\right)x_\alpha+P\frac{\partial P}{\partial{x_\alpha}} \delta_{\alpha\beta}\right]_{1\leq\alpha,\beta\leq n}}{P^{n-1}}=P^{n}.
\end{equation} 
\end{rem}

\begin{lem}\label{refine}
Without loss of generality, in the sense we clarified in Remark \ref{l1},  in \eqref{MA*} we can assume   $\lambda/2=1$. 
\end{lem}
In view of Lemma \ref{refine}, we give a refinement of Conjecture \ref{subconj}.
\begin{guess}[Refinement of Conjecture \ref{subconj}]\label{conj3}
The only projectively induced  \inv\ \KE manifolds $(M,g)$ such that $\mathrm{Ric}(g)=2g$ are open subsets of $\CP^{n_1}\times\mathellipsis\times\CP^{n_k}$, with $n_1+\mathellipsis+n_k=n$, endowed with the \K metric
$$(n_1+1) g_{FS} \oplus \mathellipsis \oplus (n_k+1) g_{FS}.$$
\end{guess}

We now sum up the results of the present section obtained so far, to arrive to Proposition \ref{equivalence} below, which plays a central role in the proof of Theorem \ref{main}, i.e. Conjecture \ref{conj3} for $n\leq 6$.

\smallskip
Every \K metric admits an infinite number of local potentials defined in the same open subset, but, in the real analytic case, there exists only one whose power expansion does not contain any term in $z$ or $\bar z$: the Diastasis function $ \dd_0(z)$.
If $(M,g)$ is a projectively induced \inv\ \K manifold and  $p$ is a fixed point of the toric action, chosen a holomorphic coordinate system $z=(z_1,\mathellipsis,z_n)$ centered at $p$, then $\dd_0(z)$ needs to read as \eqref{diastbochpi}, as seen in Lemma \ref{diastpolin}. If we in addition also assume that $(M,g)$ is a \KE manifold, then its diastasis  function $\dd_0$ can be obtained from a solution $P$ to the equation \eqref{MAeq} in the way described in Remark \ref{l1}.
 
 \begin{prop}\label{equivalence}
There exists a bijective correspondence between projectively induced  \inv\ \KE metrics defined in a neighborhood of a fixed point and the solutions of the real $n$-dimensional Monge-Ampère equation \eqref{mongeampere2} of type \eqref{type2}.
\end{prop}
\begin{proof}
As we have already explained above,  projectively induced  \inv\ \KE metrics can be obtained by finding solutions of type  \eqref{eq.Px}  to \eqref{MAeq}, and vice versa. Moreover, taking into account the replacement
 $$x_i\mapsto e^{x_i}$$
one realizes that \eqref{eq.Px} is a solution to
 \eqref{MAeq} if and only if 
\eqref{type2}
 is a solution to \eqref{mongeampere2}.
\end{proof}

Note that, in the $1$-dimensional case, \eqref{mongeampere2} reads as
$$u''=e^{-u}$$
and the unique solution of type \eqref{type2} to this equation is
\begin{equation}\label{1sol}
\log\left( 1+\frac{e^x}{2}\right)^2-x.
\end{equation}
 By substituting $x$ with $|z|^2$ in \eqref{1sol}, we get a local \K potential in the affine coordinate $z$ for the metric $2 g_{FS}$ on $\CP^1$, according to Remark \ref{uniform} and to Conjecture \ref{conj3}.

\subsection{Gradient maps of the solutions to the Monge-Ampère equation \eqref{mongeampere2}: Momentum maps and Delzant polytopes}\label{gradient}

\begin{lem}\label{polytope}
Let $u$ be a function reading as \eqref{type2}. Then the  closure $\P$ of the image $D u(\R^n)$  of the gradient map $Du$ of $u$ is the convex hull of $\I$ translated by  $-\mathbf{1}=(-1,\mathellipsis,-1)$  and, if $d=\max_{I\in\I}|I|$, 
 \begin{equation}\label{polytopeposition}
\P\subseteq\left\{(u_1,\mathellipsis,u_n)\in\R^n\,\Big|\, u_1\geq -1,\mathellipsis,u_n\geq -1,\sum_{i=1}^nu_i\leq  d-n\right\}.\end{equation}
In particular, $\P$ is a lattice polytope,  namely a polytope having all the vertices with integer coordinates.
\end{lem}
\begin{proof}
Via  a straightforward computation, we see that the values $D u(x)$ of the gradient  function are (up to translations)  convex combinations of the elements of $\I$:
\begin{equation}\label{convexcomb}
D u(x)=\frac{1}{\sum_{I\in\I} a_I e^{I\cdot x}}\sum_{I\in\I} a_I e^{I\cdot x}I-\mathbf{1}.
\end{equation}
Thus, we get that $\P+\mathbf{1}$ is the convex hull of $\I$. Moreover, since for any $I\in\I\subset \N^n$ we have that $|I|\leq d$, then for any $(u_1,\mathellipsis,u_n)\in\P+\mathbf{1}$ we have $\sum_{i=1}^nu_i\leq  d$.
\end{proof}
\begin{rem}
If $u$ is of type \eqref{type2}, then 
$$\{(-1,\mathellipsis,-1), (0,-1,\mathellipsis,-1),(-1,0,-1,\mathellipsis,-1),\mathellipsis, (-1,\mathellipsis,-1,0)\}\subset\P.$$
\end{rem}

When $u$ is a solution of \eqref{mongeampere2}, some additional geometric properties of its associated polytope $\P$ can be obtained by using the theory of toric manifolds. In order to better introduce such properties we need some extra definitions coming from the context of convex polytopes.

\begin{dfn}
Let $u:\R^n\to\R$ be a solution of type \eqref{type2} to the Monge-Ampère \eqref{mongeampere2}. The (convex) polytope $\P$ defined as the closure of $Du(\R^n)$ is going to be called the \emph{polytope associated to $u$}.
\end{dfn}

\begin{dfn}\label{delzant}
An $n$-dimensional convex polytope is  called  a \emph{Delzant} polytope if and only if all the following properties are fulfilled
\begin{description}
\item[simplicity]  at each vertex $p$, $n$ edges $l_i$ meet and
\item[rationality] $l_i=p+t v_i$, where $t\in\R^+$ and $(v_1,\mathellipsis,v_n)\in \Z^n$,
\item[smoothness] in particular $(v_1,\mathellipsis,v_n)$ is a $\Z$-basis for $\Z^n$.\end{description}\end{dfn}

\begin{dfn}\label{ref}
An $n$-dimensional  lattice polytope $\P$ containing $\mathbf{0}=(0,\mathellipsis,0)$ as an interior point is called  \emph{reflexive} if and only if 
\begin{equation}\label{A}
\P=\{ y\in\R^n\ |\ Ay\leq \mathbf{1}\},\end{equation}
where $A\in \Z^{m,n}$ and $\mathbf{1}$ is the column of length $m$ with all entries equal to $1$.
\end{dfn}
\begin{rem}\label{remarkvertex}
A reflexive polytope possesses a unique interior point with integer coordinates, which is forced to be origin in view of Definition \ref{ref}.
 \end{rem}

In order to prove that polytopes associated to solutions of type \eqref{type2} to Monge-Ampère equation \eqref{mongeampere2} are in particular Delzant (see Proposition \ref{reflexive} below), we need the following two Lemmas.

\begin{lem}\label{compact}
{A  projectively induced \inv\  \KE manifold  is an open subset of  a compact, simply connected and complete manifold $M$}.\end{lem}\begin{proof}
D. Hulin proved in \cite{hulin} that every \KE manifold \K immersed into a complex projective space can be extended to a  complete \KE manifold,  that  is also \K immersed into the same complex projective space.
Therefore, since $\lambda$ is positive by Lemma \ref{einsteinconst}, in view of the Bonnet-Myers' theorem, $M$ has to be   {compact}. Moreover, every compact \K manifold with positive definite Ricci tensor is simply connected by a well known result due to Kobayashi \cite{koricci}.
\end{proof}

\begin{lem}\label{toric}
A projectively induced  \inv\ \KE manifold $M$ is an open subset of a  toric \K manifold.
\end{lem}
\begin{proof}
In view of Lemma \ref{compact}, we can assume without loss of generality, that  $M$ is a compact, simply connected and complete  \KE manifold. Let  $z=(z_1,\mathellipsis, z_n)$ be Bochner's coordinates centered at a fixed point of the toric action and let $\dd_0(z)$ be the diastasis function.
Let $u:U\cap(\C\setminus\{0\})^n \to \R$ such that 
\begin{equation}\label{u}
u(\log{|z_1|^2},\mathellipsis, \log{|z_n|^2})=\dd_0\left(|z_1|^2,\dots,|z_n|^2\right).\end{equation}
Chosen any branch of the complex logarithm,  we  set holomorphic coordinates $w_i=\log z_i$. Hence,  the \K form reads locally as
\begin{equation}\label{moment}
\omega=\frac{i}{2} \de\bar \de u=\frac{i}{2}\sum_{k,j}\frac{\de^2 u}{\de w_k \de \bar w_j}dw_k\wedge d\bar w_j=\sum_{k,j}\frac{\de^2 u}{\de r_k \de  r_j}dr_k\wedge d \theta_j, \end{equation}
where $w_k=r_k+i \theta_k$.

Being $M$ simply connected and real analytic, each Killing vector field  $\de_{\theta_k}$ can be extended to a unique Killing vector field defined on the whole manifold $M$ (see \cite{nomizu}  theorems 1 and 2). Let $X_k$ be the global extension of $\de_{\theta_k}$.  
Therefore, since  every Killing vector field on a compact \K manifold is real holomorphic (see e.g.  \cite{Mor} prop. 9.5) and being $M$  complete, we get a holomorphic and isometric action of $\R$ on $M$ by means of the flows of the each Killing vector field $X_k$.

 Furthermore, for any $1\leq k,j\leq \dim M$, $[ X_k, X_j]$ is a Killing vector field vanishing on $U$. Since Killing vector fields (different from the identically zero vector field) vanish on totally geodesic submanifolds of real codimension at least 2 (cfr.  \cite{kob}), the commutator $[ X_k,  X_j]$ needs to vanish everywhere on $M$. Therefore, we have a transitive holomorphic and isometric action  $\mathcal{G}$ of $\R^n$ on $M$.

Let $V\subset M$ be the subset where at least a Killing vector field $X_k$ vanishes. By considering that $M\setminus V$ consists of and only of maximal dimensional orbits of the action $\mathcal{G}$, such action can be restricted to $M\setminus V$. We easily see that $\Z^n$ is the stabilizer of $\mathcal{G}$ in $(M\setminus V)\cap U$. Hence, we have a holomorphic and isometric action of the real torus  $\R^n/\Z^n$ on $M\setminus V$. Since  stabilizers of $\mathcal{G}$ in points belonging to $V$ contain $\Z^n$, this toric  action can be extended to the whole manifold $M$.

Every Killing vector field on a compact and simply connected \K manifold $M$ is Hamiltonian. Indeed, since these vector fields are also real holomorphic because $M$ is compact, they need  to be symplectic too, i.e. $i_{X_k}\omega$ is closed. Being $M$ simply connected, $H^1_{\text{dR}}(M)=0$. Therefore,  $i_{X_i}\omega$ needs to be also  exact.

 The existence of such  effective Hamiltonian action allows us to conclude that $M$ is a toric \K manifold.
\end{proof}

\begin{prop}\label{reflexive}
If a function $u$ of type \eqref{type2} is a solution to the Monge-Ampère equation \eqref{mongeampere2}, then its associated polytope  $\P$ is  Delzant and reflexive.
\end{prop} 
\begin{proof}
By considering Lemma \ref{toric}, $u$ can be seen as a  local \K potential defined on an open dense subset of an $n$-dimensional  toric \KE manifold  $M$  (see \eqref{u}). Let $\omega$ be the \K form of $M$. Since, by considering \eqref{moment}, we have
$$i_{\de_{\theta_j}}\omega=-d\left(\frac{\de u}{ \de  r_j}\right)$$
a  momentum map
$$\mu: M\to \mathfrak{t}^*\cong \R^n$$
is given  by  the gradient of $ u$ (we are denoting   the dual of the Lie algebra of $\T^n$ with $\mathfrak{t}^*$).

By the  results of T. Delzant (see e.g. \cite{del}), $\mu(M)\subset \R^n$ is forced to satisfy the  properties in the Definition \ref{delzant}.

Furthermore, since $u$ is a solution of the Monge-Ampère equation \eqref{mongeampere2}, the Ricci form $\rho$ of $M$ is equal to $2\omega$ (see also proof of Lemma \ref{einsteinconst} and Remark \ref{l1}). Then, the first Chern class $c_1(M)=\frac{1}{2\pi}[\rho]$ of $M$ is equal to $\frac{1}{\pi}[\omega]$. By taking into account also that   $\P=\mu(M)$ is  a lattice polytope (see Lemma \ref{polytope}), it follows from D. McDuff's result \cite{mcduff} that $\P$ contains  only one interior point   with integer coordinates. By \eqref{polytopeposition},  such interior point is forced to be  $\mathbf{0}$. Furthermore, in \cite{mcduff},  it is shown that the condition $c_1(M)=\frac{1}{\pi}[\omega]$ implies that the polytope $\P$ is reflexive. This has been seen by proving that  what the author called the \emph{affine distance} of the integer interior point (in our case $\mathbf{0}$) from any facet is equal to 1.
\end{proof}

Although it is widely known that the existence of  \KE metrics on toric manifolds is related to the position of the barycenter of the  image of the momentum map related to the symplectic structure (see e.g. \cite{mabuchi, bb, wangzhu}), we prove the following for the sake of self-consistency.

\begin{lem} \label{bar}
Let $u$ be a solution  of type \eqref{type2} to the Monge-Ampère equation \eqref{mongeampere2}. Then its associated polytope $\P$ has the barycenter at the origin $\mathbf{0}$.
\end{lem}
\begin{proof}
%
The $\alpha$-th component of the barycenter of $\P$ reads
\begin{multline*}
\int_{\P}u_\alpha\ du_1\mathellipsis d u_n=\int_{\R^n}\frac{\de u}{\de x_\alpha}\det D^2 u\ dx =\\
\int_{\R^n}\frac{\de u}{\de x_\alpha}e^{-u}\ dx=
-\int_{(\R^+)^{n-1}}\int_{\R^+}\frac{\de}{\de y_\alpha}\left( \frac{y_\alpha}{P(y_1,\mathellipsis,y_n)}\right)dy_\alpha d\mathbf{y},
\end{multline*}
where $y_i=e^{x_i}$ and $\mathbf{y}=(y_1,\mathellipsis,y_{\alpha -1},y_{\alpha +1},\mathellipsis, y_n )$. Since $\mathbf{0}$ is an interior point of $\P$ (see Proposition \ref{reflexive}), $$\lim_{y_\alpha\to+\infty}\frac{y_\alpha}{P(y_1,\mathellipsis,y_n)}=0,$$
 so all the previous integrals needs to vanish.\end{proof}
%


Actually we have also a converse of Lemma \ref{bar}, that holds  in the more general setting of convex bodies.
\begin{prop}[\cite{bb}]  \label{berman} If $\P$ is a convex body containing $\mathbf{0}$ in its interior, then there exists  a smooth convex function  $\phi$ satisfying the Monge-Ampère equation \eqref{mongeampere2}
 and such that the closure of the image of the gradient  map $D\phi$ of $\phi$ is $\P$ if and only if $\mathbf{0}$ is the barycenter. The solution $\phi$ is uniquely determined up to the action of the additive group $\R^n$ by translations.
\end{prop}
%
%
%
%
%

It is natural to ask the relationship between separability of solutions to the equation \eqref{mongeampere2}, i.e. solutions of type $u(x_1,\mathellipsis,x_k)+v(x_{k+1},\mathellipsis,x_{k+h})$, and decomposability of their associated polytopes, i.e. polytopes which are cartesian product of lower dimensional ones. This aspect is clarified by the following propositions.

\begin{prop}\label{product}
Let $u(x_1,\mathellipsis,x_k)$ and $v(x_{k+1},\mathellipsis,x_{k+h})$ be respectively a solution  of type \eqref{type2} to the $k$ and $h$ dimensional Monge-Ampère equation \eqref{mongeampere2}, with associated polytopes $\P$ and $\mathcal{Q}$. Then, $u+v$ is a solution of type \eqref{type2} to $k+h$ dimensional Monge-Ampère equations \eqref{mongeampere2}, whose associated polytope is $\P\times \mathcal{Q}$.
\end{prop}
\begin{proof} It follows straightforwardly.\end{proof}

Conversely, we have the following
\begin{prop}\label{product2}
If $u$ is a solution of type \eqref{type2} to the Monge-Ampère equation \eqref{mongeampere2} whose associated polytope can be decomposed as a cartesian product of a $k$-dimensional polytope $\P$ and a $h$-dimensional polytope $\mathcal{Q}$, then the $k$ and $h$ dimensional Monge-Ampère equations $\eqref{mongeampere2}$ admit solutions of type \eqref{type2} whose associated polytope are respectively $\P$ and $\mathcal{Q}$.
\end{prop}
\begin{proof}
Since $u$ is a solution of type \eqref{type2} to the $k+h$ Monge-Ampère equation \eqref{mongeampere2}, $Du(\R^{k+h})=\P\times \mathcal{Q}$ is a reflexive polytope (see Proposition \ref{reflexive}). Hence,  it is easily to see that both $\P$ and $\mathcal{Q}$ needs to be reflexive. Moreover, by considering that the origin of $\R^{k+h}$ is the  barycenter of $\P\times \mathcal{Q}$  (see Lemma \ref{bar}), the barycenter of $\P$ and $\mathcal{Q}$ is forced to be their only interior point with integer coordinates.
Therefore, in view of Proposition \ref{berman} we have convex solutions $f_1$ and $f_2$ respectively to $k$ and $h$ dimensional Monge-Ampère \eqref{mongeampere2}, whose associated polytope is respectively $\P$ and $\mathcal{Q}$. By taking into account again Proposition \ref{berman}, we have that 
$$u(x_1,\mathellipsis, x_{k+h})=f_1(x_1+c_1,\mathellipsis,x_k+c_k)+f_2(x_{k+1}+c_{k+1},\mathellipsis,x_{k+h}+c_{k+h})$$ 
with $(c_1,\mathellipsis,c_{k+h})\in\R^{k+h}$. Then $f_1(x_1+c_1,\mathellipsis,x_k+c_k)$ and $f_2(x_{k+1}+c_{k+1},\mathellipsis,x_{k+h}+c_{k+h})$ are functions of type \eqref{type2}.
\end{proof}

\subsection{Some technical results}\label{technical}
For practical reasons, in the following Lemma we refer to the Monge-Ampère \eqref{MAeq} instead of the equation \eqref{mongeampere2}, keeping always in mind the equivalence expressed by Proposition \ref{equivalence}. The following Lemma follows from Lemma 2.8 in  \cite{ny} and it will be useful for determining the shape of the  polytopes associated to the solutions of type \eqref{type2} to the Monge-Ampère equation \eqref{mongeampere2}.

\begin{lem}\label{conditions}
Let $ P$ be a solution of type \eqref{eq.Px} to the Monge-Ampère equation \eqref{MAeq}. Then the restriction $P(0,\mathellipsis,0,t,0,\mathellipsis,0)$ of $P$ to the $i$-axis is
\begin{equation}\label{P}
P(0,\mathellipsis,0,t,0,\mathellipsis,0)=\left( 1+\frac{t}{k_i}\right)^{k_i}
\end{equation}
 and the restriction $\frac{\partial P}{\partial y_j}(0,\mathellipsis,0,t,0,\mathellipsis,0)$ of $\frac{\partial P}{\partial y_j} $ to the $i$-axis, where $j\neq i$, is
\begin{equation}\label{dP}
\frac{\partial P}{\partial y_j}(0,\mathellipsis,0,t,0,\mathellipsis,0)= \left( 1+\frac{t}{k_i}\right)^{h_{ij}},
\end{equation}
for some $k_i\in\Z^+$ and $h_{ij}\in\N$. Moreover,
\begin{equation}\label{rel1}
\sum_{\alpha\neq i} h_{i\alpha }=k_i(n-2)+2
\end{equation}
and
\begin{equation}\label{rel2}
\frac{h_{ij}}{k_i}=\frac{h_{ji}}{k_j}
\end{equation}
for any $i$ and $j$.
\end{lem}
\begin{proof}
Formula \eqref{P} is formula (13) in \cite{ny} with $s=1$, where $s$ denotes the constant $\lambda /2$ with $\lambda$ the Einstein constant. Therefore, it  can be put equal to 1 in view of Lemma \ref{l1}.\\
Concerning the formula \eqref{dP}, the formula (15)  of \cite{ny}  with $s=1$ gives
$$\prod_{j\neq i} \frac{\partial P}{\partial y_j}(0,\mathellipsis,0,t,0,\mathellipsis,0)=
\prod_{\alpha=1}^R \left( 1+\frac{t}{r_i}\right)^{k_i(n-2)+2}$$
for some $R\in\Z^+$ and $r_i\in\R^+$. In the same Lemma 2.8 of  \cite{ny}, it has been proved that the only possibility is $R=1$ and $r_1=k_i$.  Hence,
$$\prod_{j\neq i} \frac{\partial P}{\partial y_j}(0,\mathellipsis,0,t,0,\mathellipsis,0)= \left( 1+\frac{t}{k_i}\right)^{k_i(n-2)+2}.$$
Then,  formulas \eqref{dP} and  \eqref{rel1} follow straightforwardly.\\
Formula \eqref{rel2} follows from Cauchy-Schwarz Lemma by  evaluating the first derivative of \eqref{dP} at the origin.
\end{proof}

\subsubsection{A geometric interpretation of the constants $k_i$ and $h_{ij}$}\label{compute}
Let $u$ be a function of type \eqref{type2}, namely a function reading as
\begin{equation}\label{U}
u(x)=\log\sum_{I\in\I} a_I e^{I\cdot x}-\sum_\alpha x_\alpha, \quad a_I=1 \text{  if  } |I|\leq 1,\end{equation}
 whose gradient image is equal to the interior of a given  polytope  $\P$. \\
Considering \eqref{convexcomb}, a direct computation shows  that the limit of $\frac{\partial u}{\partial x_j}$ for every $x_\alpha$ different from $x_i$ tending to $-\infty$, is
$$
\begin{cases}
L_i(x_i)=\frac{1}{\sum_{I\in\hat\I} a_I e^{I\cdot x}}\sum_{I\in\hat\I} a_I e^{I\cdot x}I_i-{1}, & \text{if }j=i\\
-1, & \text{otherwise}
\end{cases}$$
where $\hat\I=\{ I=(I_1,\mathellipsis,I_n) \in\I\, |\, I_\alpha=0\  \forall \alpha\neq i\}.$
Therefore the limit of the gradient of $u$, for every $x_\alpha$ different from $x_i$ tending to $-\infty$, provides a parametrization for the interior of the edge $l_i$ of $\P$ starting from $-\mathbf{1}$ and parallel to the $i$-axis. Furthermore, we have that
\begin{equation*}
\lim_{x_i\to+\infty}L_i(x_i)=\max_{I\in\hat\I} |I|-1.
\end{equation*}
Indeed, if $\hat I \in\hat \I$ is such that $|\hat I|=\max_{I\in\hat\I} |I|$, then the coefficient $a_{\hat I}$ cannot be 0 in view of the bijective correspondence between  $\I$ and integer points of $\P$ expressed by Lemma \ref{polytope}. 
By working with the  polynomial 
$$P(x)= \sum_{I\in\I} a_I x^{I},  \quad a_I=1 \text{  if  } |I|\leq 1, $$
instead of the function $u$ \eqref{U} (this choice is justified by Proposition \ref{equivalence}) we have that the degree of the restriction  to $i$-axis of $P$ is  $\max_{I\in\hat\I} |I|$, which is in turn equal to the length of the edge $l_i$ of $\P$.\\
Moreover,  by means of very similar considerations as above that we skip for the sake of brevity, we get that the degree of restriction  to the $i$-axis of the derivative of $P$ with respect to the $j$-th variable, is an integer value between 0 and the length of the intersection of $\P\cap \ell_{ij}$, where $\ell_{ij}$ denotes the straight line  parallel to the $i$-axis and passing through the point having all its coordinates equal to $-1$ except for the $j$-th one, which is equal to $0$.

\subsection{Classification of smooth reflexive polytopes: final steps of the proof of Theorem \ref{main}}\label{final}

M. Øbro developed an algorithm (see \cite{obro}) that has been used to completely classify smooth reflexive polytopes up to dimension 7. Indeed, until then  a classification only up to size 5 was known (\cite{bat1,watanabe,bat2, kn}). However, we are going to consider only polytopes up to dimension 6 because
only in this case is present a description \cite{obro2} in terms of the matrix $A$, see \eqref{A}.

Let 
\begin{equation}\label{k}
\mathbf{k}=(k_1,\mathellipsis,k_n)^T\end{equation}
 and
\begin{equation}\label{H} 
\mathbf{H}_{ij}=\begin{cases}
h_{ij} & \text{if }i\neq j\\
0 & \text{if } i=j
\end{cases}\end{equation}
where $k_i$ and $h_{ij}$ are those defined in Lemma \ref{conditions}. 

By Propositions \ref{reflexive}, \ref{bar} and \ref{product}, we can consider only
Delzant reflexive polytopes with barycenter at $\mathbf{0}$ that cannot be decomposed as a cartesian product of lower dimensional polytopes. 

\smallskip
In the subsequent subsections we will use some tables containing the matrix $A$, defined  by \eqref{A},
and the vector $\mathbf{k}$ and the matrix $\mathbf{H}$, defined respectively by \eqref{k} and \eqref{H}, that we have computed by considering what seen in Section \ref{compute}. Notice that we consider only the case where each entry of $\mathbf{H}$ attains the maximum value predicted in the aforementioned section, because we are going to realize that only in this case the condition \eqref{rel1} is satisfied.

\subsubsection{Simplex}\label{simplex}
The $n$-simplex with the following $n+1$ vertices 
$$(-1,\mathellipsis,-1),\, (n,-1,\mathellipsis,-1),\, (-1,n,-1\mathellipsis),\,\mathellipsis, \, (-1,\mathellipsis,-1,n)$$
is Delzant and reflexive and, as such, there exists a unique convex solution of the Monge-Ampère \eqref{mongeampere2} associated to it.  The aforementioned simplex is described by Table \ref{table1}.

\begin{longtable}{ccc}
\toprule
$A$ & $\mathbf{k}$ &  $\mathbf{H}$\\
\toprule
$\left(
\begin{array}{c}
 -\mathrm{Id}_{n\times n}\\
 \mathbf{1}  \\
\end{array}
\right)$ & $
\left(
\begin{array}{c}
 n+1\\
 \vdots \\
n+1\\
\end{array}
\right)$ & $\mathbf{H}_{ij}=n-n\delta_{ij}$\\\arrayrulecolor{white}
 \midrule
\caption{$n$-dimensional simplex with barycenter at the origin.}
\label{table1}
\end{longtable}
In Table \ref{table1}, $\mathrm{Id}_{n\times n}$ is the $n$-dimensional identity matrix and $\mathbf{1}$ is the row whose entries are 1.
We note that the values contained in the table do not contradict \eqref{rel1} and \eqref{rel2}. Indeed, in this case we have the unique solution
\begin{equation}\label{cpn}
u(x_1,\mathellipsis,x_n)=\log\left(1+\sum_{i=1}^n \frac{e^{x_i}}{n+1}\right)^{n+1}-\sum_{i=1}^n x_i,
\end{equation}
that, in view of Proposition \ref{equivalence}, is the solution associated to 
\begin{equation}\label{cpn2}
\left(\CP^n, (n+1)  g_{FS}\right),\end{equation}
in accordance  with Conjecture \ref{conj3}.

\subsubsection{1-dimensional case}\label{1dim}

As we have already seen in the very end of Section \ref{reduction}, the only solution in this case is  \eqref{1sol}, that leads to \eqref{cpn2} for $n=1$. In particular, we notice that the associated polytope is  the segment from $-1$ to $1$, namely a $1$-dimensional simplex, according to Section \ref{simplex}.

\subsubsection{2-dimensional case}\label{2dim}
As we have already seen in Section \ref{simplex}, in this case we have the solution associated to the 2-simplex, namely \eqref{cpn} for $n=2$. Furthermore, in view of the Proposition \ref{product} and Section \ref{1dim}, we have also the solution
$$\log\left( 1 +\frac{x_1}{2}\right)^2+ \log\left( 1 +\frac{x_2}{2}\right)^2-x_1-x_2,$$
whose associated polytope is the square with vertices $(-1,-1)$, $(-1,1)$,  $(1,-1)$ and $(1,1)$, namely the cartesian product of the segment $\{(t,-1)\ |\ -1\leq t \leq 1\}$ and the segment  $\{(-1,t)\ |\ -1\leq t \leq 1\}$. This is the only reflexive Delzant polytope that can be decomposed as a cartesian product of 1-dimensional ones. In view of Proposition \ref{equivalence}, this solution leads to $$(\CP^1\times\CP^1,\ 2g_{FS}\oplus 2g_{FS}),$$
cfr. Conjecture \ref{conj3}.
There exists another reflexive Delzant polytope with barycenter at the origin, namely the one given by Table \ref{table2}.
\begin{longtable}{ccc}
\arrayrulecolor{black}
\toprule
$A$ & $\mathbf{k}$ &  $\mathbf{H}$\\
\toprule
$ \begin{psmallmatrix}
 -1 & 0 \\
 0 & -1 \\
 1 & -1 \\
 -1 & 1 \\
 1 & 0 \\
 0 & 1 \\
\end{psmallmatrix}
$ & $ \begin{psmallmatrix}
 1 \\
 1 \\
\end{psmallmatrix}
$ & $ \begin{psmallmatrix}
 0 & 2 \\
 2 & 0 \\
\end{psmallmatrix}
$\\\arrayrulecolor{white}
 \midrule
\caption{Undecomposable 2-dimensional smooth reflexive polytopes with barycenter at the origin (2-simplex  excluded).}
\label{table2}
\end{longtable}
The polytope described by Table \ref{table2} is the hexagon $\mathcal{E}$ with vertices $(-1,-1)$, $(-1,0)$, $(0,-1)$, $(0,1)$, $(1,0)$, $ (1,1)$. It is easy to realize that such polytope satisfies conditions \eqref{rel1} and \eqref{rel2}. 
The most general function $u$ of type \eqref{type2} such that the closure of $Du(\R^2)$ is equal to $\mathcal{E}$ reads as
$$\log\left(1+e^{x_1}+e^{x_2}+a_{(1,1)} e^{x_1+x_2} +a_{(2,0)} e^{2x_1+x_2}+a_{(0,2)} e^{x_1+2x_2}+a_{(2,2)} e^{2x_1+2x_2}\right)-x_1-x_2.$$
If such $u$ is a solution of the Monge-Ampère \eqref{mongeampere2}, then it needs to satisfy \eqref{P} and \eqref{dP}. Hence
$$u(x_1,x_2)=\log\left(1+e^{x_1}+e^{x_2}+2 e^{x_1+x_2} + e^{2x_1+x_2}+ e^{x_1+2x_2}+a_{(2,2)} e^{2x_1+2x_2}\right)-x_1-x_2.$$
We can easily compute that
$$\lim_{x_1\to-\infty}\lim_{x_2\to -\infty} \frac{\partial^2}{\partial x_1\partial x_2}\left(e^u \det D^2 u\right)\neq 0$$
independently of $a_{(2,2)}$. Therefore there is no  $a_{(2,2)}\in \R$ for which  $u$ is a solution of \eqref{mongeampere2}.

\subsubsection{3-dimensional case}\label{3dim}
As said in Section \ref{simplex}, we have the solution associated to the 3-simplex, namely \eqref{cpn} for $n=3$.
 Furthermore, in view of the Proposition \ref{product} and Section \ref{1dim}, we have also the solutions
$$\log\left( 1 +\frac{x_1}{2}\right)^2+ \log\left( 1 +\frac{x_2}{2}\right)^2+\left( 1 +\frac{x_3}{2}\right)^2-x_1-x_2-x_3$$
and, up to variables renaming,
$$\log\left( 1 +\frac{x_1+x_2}{3}\right)^3+ \log\left( 1 +\frac{x_3}{2}\right)^2-x_1-x_2-x_3,$$
which, in view of Proposition \ref{equivalence}, lead respectively to 
$$(\CP^1\times\CP^1\times \CP^1,\ 2g_{FS}\oplus 2g_{FS}\oplus 2g_{FS})$$
and 
$$(\CP^2\times\CP^1,\ 3g_{FS}\oplus 2g_{FS}),$$ cfr. Conjecture \ref{conj3}.
We can easily see that the associated polytopes are respectively, the cube with vertices 
$$(-1,-1,-1), (1,-1,-1), (-1,1,-1), (-1,-1,1), (1,1,-1), (1,-1,1), (-1,1,1), (1,1,1),$$
namely the cartesian product of the three segments
$\{(-1,-1,t)\ |\ -1\leq t \leq 1\}$, $\{(-1,t,-1)\ |\ -1\leq t \leq 1\}$ and $\{(t,-1,-1)\ |\ -1\leq t \leq 1\}$,\\
and the prism with vertices 
$$(-1,-1,-1), (-1,1,-1), (1,-1,-1), (-1,-1,1), (-1,1,1), (1,-1,1)$$
namely the cartesian product of the 2-simplex whose vertices are $(-1,-1,-1), (-1,1,-1), (1,-1,-1)$ and the segment $\{(-1,-1,t)\ |\ -1\leq t \leq 1\}$. By taking into account Proposition \ref{product2}, we have no more decomposable polytopes to take into account. Indeed, even if there exists another decomposable reflexive Delzant polytope with barycenter at the origin, namely the prism with vertices
$$(-1,-1,-1), (-1,0,-1), (0,-1,-1), (0,1,-1), (1,0,-1),  (1,1,-1),$$
$$(-1,-1,1), (-1,0,1), (0,-1,1), (0,1,1), (1,0,1),  (1,1,1),$$
we can see directly that it is a cartesian product of an hexagon and a segment.
In view of Proposition \ref{product2}, there are no solutions of type \eqref{type2}  associated to such polytope, since there are no 2-dimensional solutions of type \eqref{type2} associated to the hexagon.

Finally, there is another undecomposable reflexive Delzant polytope with barycenter at the origin, namely the one given by Table \ref{table3}.
\begin{longtable}{ccc}
\arrayrulecolor{black}
\toprule
$A$ & $\mathbf{k}$ &  $\mathbf{H}$\\
\toprule
$ \begin{psmallmatrix}
 -1 & 0 & 0 \\
 0 & -1 & 0 \\
 0 & 0 & -1 \\
 1 & 0 & -1 \\
 0 & 0 & 1 \\
 0 & 1 & 1 \\
\end{psmallmatrix}
$ & $\begin{psmallmatrix}
 1 \\
 3 \\
 2 \\
\end{psmallmatrix}
$ & $\begin{psmallmatrix}
 0 & 1 & 2 \\
 3 & 0 & 2 \\
 2 & 2 & 0 \\
\end{psmallmatrix}
$\\\arrayrulecolor{white}
 \midrule
\caption{Undecomposable 3-dimensional smooth reflexive polytopes with barycenter at the origin (3-simplex  excluded).}
\label{table3}
\end{longtable}
Note that the values contained in Table \ref{table3}  do not satisfy the condition \eqref{rel2}. Therefore we cannot have solutions of type \eqref{type2} related to such polytope.

\subsubsection{4-dimensional case}\label{4dim}
As said in Section \ref{simplex}, we have the solution associated to the 4-simplex, namely \eqref{cpn} for $n=4$.
Moreover, in view of the Proposition \ref{product2}, Proposition \ref{product} and taking into account the results of Sections \ref{1dim} - \ref{3dim}, the only solutions whose associated polytope can be decomposed as a cartesian product of lower dimensional polytopes are
$$\log\left( 1 +\frac{x_1}{2}\right)^2+ \log\left( 1 +\frac{x_2}{2}\right)^2+\left( 1 +\frac{x_3}{2}\right)^2+\log\left( 1 +\frac{x_4}{2}\right)^2-x_1-x_2-x_3-x_4,$$
$$\log\left( 1 +\frac{x_1+x_2}{3}\right)^3+\left( 1 +\frac{x_3}{2}\right)^2+\log\left( 1 +\frac{x_4}{2}\right)^2-x_1-x_2-x_3-x_4,$$
$$\log\left( 1 +\frac{x_1+x_2}{3}\right)^3+\log\left( 1 +\frac{x_3+x_4}{3}\right)^3-x_1-x_2-x_3-x_4,$$
$$\log\left( 1 +\frac{x_1+x_2+x_3}{4}\right)^4+\log\left( 1 +\frac{x_4}{2}\right)^2-x_1-x_2-x_3-x_4.$$
In view of Proposition \ref{equivalence}, these solutions respectively lead to
$$(\CP^1\times\CP^1\times \CP^1\times \CP^1,\ 2g_{FS}\oplus 2g_{FS}\oplus 2g_{FS}\oplus 2g_{FS}),$$
$$(\CP^2\times\CP^1\times \CP^1,\ 3g_{FS}\oplus 2g_{FS}\oplus 2g_{FS}),$$
$$(\CP^2\times\CP^2,\ 3g_{FS}\oplus 3g_{FS}),$$
$$(\CP^3\times\CP^1,\ 4g_{FS}\oplus 2g_{FS}),$$ cfr. Conjecture \ref{conj3}.
Beside the 4-simplex, there exists also three further undecomposable reflexive Delzant polytope with barycenter at the origin, namely the ones given by Table \ref{table4}.
Only the first polytope in such table satisfies both conditions \eqref{rel1} and \eqref{rel2}. Nevertheless, if we assume the existence of a solution of type \eqref{type2} to the Monge-Ampère equation \eqref{mongeampere2} associated to such polytope, we get a contradiction. Indeed,
by taking into account \eqref{P} and \eqref{dP}, we obtain after long computations that
$$\lim_{x_1\to-\infty}\lim_{x_2\to -\infty}\lim_{x_3\to -\infty}\lim_{x_4\to -\infty} \frac{\partial^2}{\partial x_1\partial x_3}\left(e^u \det D^2 u\right)>0.$$

\subsubsection{5-dimensional case}\label{5dim}
As already seen in  Section
 \ref{simplex}, we have the solution associated to the 5-simplex, namely \eqref{cpn} for $n=5$. 
 Moreover, in view of Proposition \ref{product} and \ref{product2},  we can obtain all the solutions associated to  decomposable polytopes (as a cartesian product of lower dimensional ones) by taking into account the results of Sections \ref{1dim}-\ref{4dim}.
Beside the 5-simplex, there are also  seven further undecomposable reflexive Delzant 5-polytope with barycenter at the origin (see Table \ref{table5}). Nevertheless, none of them satisfies the condition \eqref{rel2}.

\subsubsection{6-dimensional case}\label{6dim}
As already seen in  Section
 \ref{simplex}, we have the solution associated to the 6-simplex, namely \eqref{cpn} for $n=6$. 
 Moreover, in view of Proposition \ref{product} and \ref{product2},  we can obtain all the solutions associated to  decomposable polytopes (as a cartesian product of lower dimensional ones) by taking into account the results of Sections \ref{1dim}-\ref{5dim}.
Beside the 6-simplex, there are also  twelve further undecomposable reflexive Delzant 6-polytope with barycenter at the origin, but only one (the first polytope in the Table \ref{table6}) satisfies both the condition  \eqref{rel1}  and \eqref{rel2}. Nevertheless, if we assume the existence of a solution of type \eqref{type2} to the Monge-Ampère equation \eqref{mongeampere2} associated to such polytope, we get a contradiction. Indeed,
by taking into account \eqref{P} and \eqref{dP}, we can consider a linear system
$$\lim_{x_1\to-\infty}\lim_{x_2\to -\infty}\lim_{x_3\to -\infty}\lim_{x_4\to -\infty}\lim_{x_5\to -\infty}\lim_{x_6\to -\infty} \frac{\partial^2}{\partial x_\alpha \partial x_\beta}\left(e^u \det D^2 u\right)=0$$
where $\alpha,\beta=1,\mathellipsis,6$,
having the  third degree coefficients of the polynomial $P=e^{u+\sum_{i=1}^6 y_i}|_{y_i=\log x_i}$ as variables. By considering that any coefficient of $P$ cannot be negative, we obtain that such system admits a unique solution. Thus, this result puts us in the position to get, after long computations,
$$\lim_{x_1\to-\infty}\lim_{x_2\to -\infty}\lim_{x_3\to -\infty}\lim_{x_4\to -\infty}\lim_{x_5\to -\infty}\lim_{x_6\to -\infty} \frac{\partial^3}{\partial x_1^2 \partial x_2}\left(e^u \det D^2 u\right)>0,$$
that clearly contradicts \eqref {mongeampere2}.

\newpage
\begin{table}
\parbox{ .45\linewidth}{
\centering
\begin{longtable}{ccc}
\arrayrulecolor{black}
\toprule
$A$ & $\mathbf{k}$ &  $\mathbf{H}$\\
\toprule
$\begin{psmallmatrix}
 -1 & 0 & 0 & 0 \\
 0 & -1 & 0 & 0 \\
 0 & 0 & -1 & 0 \\
 0 & 0 & 0 & -1 \\
 1 & 1 & -1 & -1 \\
 -1 & -1 & 1 & 1 \\
 1 & 0 & 0 & 0 \\
 0 & 1 & 0 & 0 \\
 0 & 0 & 1 & 0 \\
 0 & 0 & 0 & 1 \\
\end{psmallmatrix}
$ & $\begin{psmallmatrix}
 1 \\
 1 \\
 1 \\
 1 \\
\end{psmallmatrix}
 $ & $\begin{psmallmatrix}
 0 & 0 & 2 & 2 \\
 0 & 0 & 2 & 2 \\
 2 & 2 & 0 & 0 \\
 2 & 2 & 0 & 0 \\
\end{psmallmatrix}
 $\\\arrayrulecolor{white}
 \midrule
$\begin{psmallmatrix}
 -1 & 0 & 0 & 0 \\
 0 & -1 & 0 & 0 \\
 0 & 0 & -1 & 0 \\
 0 & 0 & 0 & -1 \\
 1 & 1 & -1 & -1 \\
 0 & -1 & 1 & 0 \\
 -1 & 0 & 0 & 1 \\
 1 & 1 & 0 & 0 \\
 0 & 0 & 1 & 1 \\
\end{psmallmatrix}
 $ & $\begin{psmallmatrix}
 1 \\
 1 \\
 1 \\
 1 \\
\end{psmallmatrix}
 $ & $\begin{psmallmatrix}
 0 & 0 & 2 & 2 \\
 0 & 0 & 2 & 2 \\
 1 & 2 & 0 & 1 \\
 2 & 1 & 1 & 0 \\
\end{psmallmatrix}
 $\\\arrayrulecolor{white}
 \midrule
$\begin{psmallmatrix}
 -1 & 0 & 0 & 0 \\
 0 & -1 & 0 & 0 \\
 0 & 0 & -1 & 0 \\
 0 & 0 & 0 & -1 \\
 1 & 0 & 0 & -1 \\
 0 & 1 & 0 & -1 \\
 0 & -1 & 0 & 1 \\
 0 & 1 & 0 & 0 \\
 0 & 0 & 0 & 1 \\
 0 & 0 & 1 & 1 \\
\end{psmallmatrix}
$ & $\begin{psmallmatrix}
 1 \\
 1 \\
 3 \\
 1 \\
\end{psmallmatrix}
$ & $\begin{psmallmatrix}
 0 & 1 & 1 & 2 \\
 1 & 0 & 1 & 2 \\
 3 & 3 & 0 & 2 \\
 1 & 2 & 1 & 0 \\
\end{psmallmatrix}
$\\\arrayrulecolor{white}
 \midrule
\caption{Undecomposable 4-dimensional smooth reflexive polytopes with barycenter at the origin (4-symplex  excluded).}
\label{table4}
\end{longtable}}
\end{table}

\begin{table}
\parbox{ .45\linewidth}{
\centering
\begin{tabular}{ccc}
\arrayrulecolor{black}
 \toprule
$A$ & $\mathbf{k}$ &  $\mathbf{H}$\\
\toprule
$ \begin{psmallmatrix}
 -1 & 0 & 0 & 0 & 0 \\
 0 & -1 & 0 & 0 & 0 \\
 0 & 0 & -1 & 0 & 0 \\
 0 & 0 & 0 & -1 & 0 \\
 0 & 0 & 0 & 0 & -1 \\
 1 & 1 & 0 & 0 & -2 \\
 0 & 0 & 0 & 0 & 1 \\
 0 & 0 & 1 & 1 & 2 \\
\end{psmallmatrix}$ & $
\begin{psmallmatrix}
 1 \\
 1 \\
 5 \\
 5 \\
 2 \\
\end{psmallmatrix}$ & $\begin{psmallmatrix}
 0 & 0 & 1 & 1 & 3 \\
 0 & 0 & 1 & 1 & 3 \\
 5 & 5 & 0 & 4 & 3 \\
 5 & 5 & 4 & 0 & 3 \\
 2 & 2 & 2 & 2 & 0 \\
\end{psmallmatrix}$\\
\arrayrulecolor{white}
 \midrule
$\begin{psmallmatrix}
 -1 & 0 & 0 & 0 & 0 \\
 0 & -1 & 0 & 0 & 0 \\
 0 & 0 & -1 & 0 & 0 \\
 0 & 0 & 0 & -1 & 0 \\
 0 & 0 & 0 & 0 & -1 \\
 1 & 1 & 0 & -1 & -1 \\
 -1 & -1 & 0 & 1 & 1 \\
 1 & 0 & 0 & 0 & 0 \\
 0 & 1 & 0 & 0 & 0 \\
 0 & 0 & 1 & 1 & 1 \\
\end{psmallmatrix}$ & $\begin{psmallmatrix}
 1 \\
 1 \\
 4 \\
 1 \\
 1 \\
\end{psmallmatrix}$ & $\begin{psmallmatrix}
 0 & 0 & 1 & 2 & 2 \\
 0 & 0 & 1 & 2 & 2 \\
 4 & 4 & 0 & 3 & 3 \\
 2 & 2 & 1 & 0 & 0 \\
 2 & 2 & 1 & 0 & 0 \\
\end{psmallmatrix}$\\
\arrayrulecolor{white}
 \midrule
$\begin{psmallmatrix}
 -1 & 0 & 0 & 0 & 0 \\
 0 & -1 & 0 & 0 & 0 \\
 0 & 0 & -1 & 0 & 0 \\
 0 & 0 & 0 & -1 & 0 \\
 0 & 0 & 0 & 0 & -1 \\
 1 & 0 & 0 & 0 & -1 \\
 0 & 1 & 0 & 0 & -1 \\
 0 & -1 & 0 & 0 & 1 \\
 0 & 1 & 0 & 0 & 0 \\
 0 & 0 & 0 & 0 & 1 \\
 0 & -1 & 1 & 0 & 1 \\
 0 & 1 & 0 & 1 & 0 \\
\end{psmallmatrix}$ & $\begin{psmallmatrix}
 1 \\
 1 \\
 2 \\
 3 \\
 1 \\
\end{psmallmatrix}$ & $\begin{psmallmatrix}
 0 & 1 & 1 & 1 & 2 \\
 1 & 0 & 1 & 1 & 2 \\
 2 & 3 & 0 & 2 & 1 \\
 3 & 2 & 3 & 0 & 3 \\
 1 & 2 & 1 & 1 & 0 \\
\end{psmallmatrix}$\\
\arrayrulecolor{white}
 \midrule
\end{tabular}}
\hfill
\parbox{ .45\linewidth}{
\centering
\begin{tabular}{ccc}
\arrayrulecolor{black}
\toprule
$A$ & $\mathbf{k}$ &  $\mathbf{H}$\\
\toprule
$\begin{psmallmatrix}
 -1 & 0 & 0 & 0 & 0 \\
 0 & -1 & 0 & 0 & 0 \\
 0 & 0 & -1 & 0 & 0 \\
 0 & 0 & 0 & -1 & 0 \\
 0 & 0 & 0 & 0 & -1 \\
 1 & 0 & 0 & 0 & -1 \\
 0 & 1 & 0 & 0 & -1 \\
 0 & 0 & 0 & 0 & 1 \\
 0 & 0 & 1 & 0 & 1 \\
 0 & 0 & 0 & 1 & 1 \\
\end{psmallmatrix}$ & $\begin{psmallmatrix}
 1 \\
 1 \\
 3 \\
 3 \\
 2 \\
\end{psmallmatrix}$ & $\begin{psmallmatrix}
 0 & 1 & 1 & 1 & 2 \\
 1 & 0 & 1 & 1 & 2 \\
 3 & 3 & 0 & 3 & 2 \\
 3 & 3 & 3 & 0 & 2 \\
 2 & 2 & 2 & 2 & 0 \\
\end{psmallmatrix}$\\\arrayrulecolor{white}
 \midrule
$\begin{psmallmatrix}
 -1 & 0 & 0 & 0 & 0 \\
 0 & -1 & 0 & 0 & 0 \\
 0 & 0 & -1 & 0 & 0 \\
 0 & 0 & 0 & -1 & 0 \\
 0 & 0 & 0 & 0 & -1 \\
 1 & 0 & 0 & 0 & -1 \\
 0 & 1 & 0 & -1 & 0 \\
 0 & 0 & 0 & 1 & 1 \\
 0 & 0 & 1 & 1 & 1 \\
\end{psmallmatrix}$ & $\begin{psmallmatrix}
 1 \\
 1 \\
 4 \\
 3 \\
 3 \\
\end{psmallmatrix}$ & $\begin{psmallmatrix}
 0 & 1 & 1 & 1 & 2 \\
 1 & 0 & 1 & 2 & 1 \\
 4 & 4 & 0 & 3 & 3 \\
 3 & 3 & 3 & 0 & 2 \\
 3 & 3 & 3 & 2 & 0 \\
\end{psmallmatrix}$\\\arrayrulecolor{white}
 \midrule
$ \begin{psmallmatrix}
 -1 & 0 & 0 & 0 & 0 \\
 0 & -1 & 0 & 0 & 0 \\
 0 & 0 & -1 & 0 & 0 \\
 0 & 0 & 0 & -1 & 0 \\
 0 & 0 & 0 & 0 & -1 \\
 1 & 1 & 0 & 0 & -1 \\
 0 & 0 & 0 & 0 & 1 \\
 0 & 0 & 1 & 1 & 1 \\
\end{psmallmatrix}$  & $\begin{psmallmatrix}
 2 \\
 2 \\
 4 \\
 4 \\
 2 \\
\end{psmallmatrix}$ & $\begin{psmallmatrix}
 0 & 1 & 2 & 2 & 3 \\
 1 & 0 & 2 & 2 & 3 \\
 4 & 4 & 0 & 3 & 3 \\
 4 & 4 & 3 & 0 & 3 \\
 2 & 2 & 2 & 2 & 0 \\
\end{psmallmatrix}$
\end{tabular}}
\caption{Undecomposable 5-dimensional smooth reflexive polytopes with barycenter at the origin (5-symplex  excluded).}
\label{table5}
\end{table}

\begin{table}
\parbox{ .45\linewidth}{
\centering
\begin{tabular}{ccc}
\arrayrulecolor{black}
\toprule
$A$ & $\mathbf{k}$ &  $\mathbf{H}$\\
\toprule
$ \begin{psmallmatrix}
 -1 & 0 & 0 & 0 & 0 & 0 \\
 0 & -1 & 0 & 0 & 0 & 0 \\
 0 & 0 & -1 & 0 & 0 & 0 \\
 0 & 0 & 0 & -1 & 0 & 0 \\
 0 & 0 & 0 & 0 & -1 & 0 \\
 0 & 0 & 0 & 0 & 0 & -1 \\
 1 & 1 & 1 & -1 & -1 & -1 \\
 -1 & -1 & -1 & 1 & 1 & 1 \\
 1 & 0 & 0 & 0 & 0 & 0 \\
 0 & 1 & 0 & 0 & 0 & 0 \\
 0 & 0 & 1 & 0 & 0 & 0 \\
 0 & 0 & 0 & 1 & 0 & 0 \\
 0 & 0 & 0 & 0 & 1 & 0 \\
 0 & 0 & 0 & 0 & 0 & 1 \\
\end{psmallmatrix}$ & $\begin{psmallmatrix}
 1 \\
 1 \\
 1 \\
 1 \\
 1 \\
 1 \\
\end{psmallmatrix}
$ & $\begin{psmallmatrix}
 0 & 0 & 0 & 2 & 2 & 2 \\
 0 & 0 & 0 & 2 & 2 & 2 \\
 0 & 0 & 0 & 2 & 2 & 2 \\
 2 & 2 & 2 & 0 & 0 & 0 \\
 2 & 2 & 2 & 0 & 0 & 0 \\
 2 & 2 & 2 & 0 & 0 & 0 \\
\end{psmallmatrix}
$\\\arrayrulecolor{white}
 \midrule
$\begin{psmallmatrix}
 -1 & 0 & 0 & 0 & 0 & 0 \\
 0 & -1 & 0 & 0 & 0 & 0 \\
 0 & 0 & -1 & 0 & 0 & 0 \\
 0 & 0 & 0 & -1 & 0 & 0 \\
 0 & 0 & 0 & 0 & -1 & 0 \\
 0 & 0 & 0 & 0 & 0 & -1 \\
 1 & 1 & 1 & -1 & -1 & -1 \\
 0 & -1 & -1 & 1 & 1 & 0 \\
 0 & 1 & 1 & 0 & -1 & 0 \\
 -1 & 0 & -1 & 1 & 1 & 1 \\
 1 & 0 & 1 & 0 & 0 & 0 \\
 0 & 0 & 0 & 0 & 1 & 1 \\
\end{psmallmatrix}
$ & $\begin{psmallmatrix}
 1 \\
 1 \\
 1 \\
 1 \\
 1 \\
 2 \\
\end{psmallmatrix}
$ & $\begin{psmallmatrix}
 0 & 0 & 0 & 2 & 2 & 2 \\
 0 & 0 & 0 & 2 & 2 & 2 \\
 0 & 0 & 0 & 2 & 2 & 2 \\
 1 & 2 & 2 & 0 & 0 & 1 \\
 1 & 2 & 2 & 0 & 0 & 1 \\
 3 & 2 & 3 & 1 & 1 & 0 \\
\end{psmallmatrix}
$\\\arrayrulecolor{white}
 \midrule
$\begin{psmallmatrix}
 -1 & 0 & 0 & 0 & 0 & 0 \\
 0 & -1 & 0 & 0 & 0 & 0 \\
 0 & 0 & -1 & 0 & 0 & 0 \\
 0 & 0 & 0 & -1 & 0 & 0 \\
 0 & 0 & 0 & 0 & -1 & 0 \\
 0 & 0 & 0 & 0 & 0 & -1 \\
 1 & 1 & 1 & -1 & -1 & -1 \\
 0 & 0 & -1 & 1 & 0 & 0 \\
 0 & -1 & 0 & 0 & 1 & 0 \\
 -1 & 0 & 0 & 0 & 0 & 1 \\
 1 & 1 & 1 & 0 & 0 & 0 \\
 0 & 0 & 0 & 1 & 1 & 1 \\
\end{psmallmatrix}
$ & $\begin{psmallmatrix}
 1 \\
 1 \\
 1 \\
 1 \\
 1 \\
 1 \\
\end{psmallmatrix}
$ & $\begin{psmallmatrix}
 0 & 0 & 0 & 2 & 2 & 2 \\
 0 & 0 & 0 & 2 & 2 & 2 \\
 0 & 0 & 0 & 2 & 2 & 2 \\
 1 & 1 & 2 & 0 & 1 & 1 \\
 1 & 2 & 1 & 1 & 0 & 1 \\
 2 & 1 & 1 & 1 & 1 & 0 \\
\end{psmallmatrix}
$\\\arrayrulecolor{white}
 \midrule
$\begin{psmallmatrix}
 -1 & 0 & 0 & 0 & 0 & 0 \\
 0 & -1 & 0 & 0 & 0 & 0 \\
 0 & 0 & -1 & 0 & 0 & 0 \\
 0 & 0 & 0 & -1 & 0 & 0 \\
 0 & 0 & 0 & 0 & -1 & 0 \\
 0 & 0 & 0 & 0 & 0 & -1 \\
 1 & 1 & 0 & 0 & 0 & -2 \\
 0 & 0 & 1 & 0 & 0 & -1 \\
 0 & 0 & -1 & 0 & 0 & 1 \\
 0 & 0 & 1 & 0 & 0 & 0 \\
 0 & 0 & 0 & 0 & 0 & 1 \\
 0 & 0 & 0 & 1 & 1 & 2 \\
\end{psmallmatrix}
$ & $\begin{psmallmatrix}
 1 \\
 1 \\
 1 \\
 5 \\
 5 \\
 1 \\
\end{psmallmatrix}
$ & $\begin{psmallmatrix}
 0 & 0 & 1 & 1 & 1 & 3 \\
 0 & 0 & 1 & 1 & 1 & 3 \\
 1 & 1 & 0 & 1 & 1 & 2 \\
 5 & 5 & 5 & 0 & 4 & 3 \\
 5 & 5 & 5 & 4 & 0 & 3 \\
 1 & 1 & 2 & 1 & 1 & 0 \\
\end{psmallmatrix}
$\\\arrayrulecolor{white}
 \midrule
$\begin{psmallmatrix}
 -1 & 0 & 0 & 0 & 0 & 0 \\
 0 & -1 & 0 & 0 & 0 & 0 \\
 0 & 0 & -1 & 0 & 0 & 0 \\
 0 & 0 & 0 & -1 & 0 & 0 \\
 0 & 0 & 0 & 0 & -1 & 0 \\
 0 & 0 & 0 & 0 & 0 & -1 \\
 1 & 1 & 0 & 0 & -1 & -1 \\
 0 & -1 & 1 & -1 & 1 & 0 \\
 -1 & 0 & -1 & 1 & 0 & 1 \\
 1 & 0 & 1 & 0 & 0 & 0 \\
 0 & 1 & 0 & 1 & 0 & 0 \\
 0 & 0 & 0 & 0 & 1 & 1 \\
\end{psmallmatrix}
$ & $\begin{psmallmatrix}
 1 \\
 1 \\
 1 \\
 1 \\
 1 \\
 1 \\
\end{psmallmatrix}
$ & $\begin{psmallmatrix}
 0 & 0 & 1 & 1 & 2 & 2 \\
 0 & 0 & 1 & 1 & 2 & 2 \\
 1 & 2 & 0 & 2 & 0 & 1 \\
 2 & 1 & 2 & 0 & 1 & 0 \\
 1 & 2 & 0 & 2 & 0 & 1 \\
 2 & 1 & 2 & 0 & 1 & 0 \\
\end{psmallmatrix}
$\\\arrayrulecolor{white}
 \midrule
$\begin{psmallmatrix}
 -1 & 0 & 0 & 0 & 0 & 0 \\
 0 & -1 & 0 & 0 & 0 & 0 \\
 0 & 0 & -1 & 0 & 0 & 0 \\
 0 & 0 & 0 & -1 & 0 & 0 \\
 0 & 0 & 0 & 0 & -1 & 0 \\
 0 & 0 & 0 & 0 & 0 & -1 \\
 1 & 1 & 0 & 0 & -1 & -1 \\
 0 & 0 & 1 & 0 & 0 & -1 \\
 -1 & -1 & 0 & 0 & 1 & 1 \\
 1 & 0 & 0 & 0 & 0 & 0 \\
 0 & 1 & 0 & 0 & 0 & 0 \\
 0 & 0 & 0 & 0 & 1 & 0 \\
 0 & 0 & 0 & 0 & 0 & 1 \\
 0 & 0 & 0 & 1 & 0 & 1 \\
\end{psmallmatrix}
$ & $\begin{psmallmatrix}
 1 \\
 1 \\
 1 \\
 3 \\
 1 \\
 1 \\
\end{psmallmatrix}
$ & $\begin{psmallmatrix}
 0 & 0 & 1 & 1 & 2 & 2 \\
 0 & 0 & 1 & 1 & 2 & 2 \\
 1 & 1 & 0 & 1 & 1 & 2 \\
 3 & 3 & 3 & 0 & 3 & 2 \\
 2 & 2 & 1 & 1 & 0 & 0 \\
 2 & 2 & 1 & 1 & 0 & 0 \\
\end{psmallmatrix}
$\\\arrayrulecolor{white}
 \midrule
\end{tabular}}
\hfill
\parbox{ .45\linewidth}{
\centering
\begin{tabular}{ccc}
\arrayrulecolor{black}
\toprule
$A$ & $\mathbf{k}$ &  $\mathbf{H}$\\
\toprule
$\begin{psmallmatrix}
 -1 & 0 & 0 & 0 & 0 & 0 \\
 0 & -1 & 0 & 0 & 0 & 0 \\
 0 & 0 & -1 & 0 & 0 & 0 \\
 0 & 0 & 0 & -1 & 0 & 0 \\
 0 & 0 & 0 & 0 & -1 & 0 \\
 0 & 0 & 0 & 0 & 0 & -1 \\
 1 & 1 & 0 & 0 & -1 & -1 \\
 -1 & -1 & 0 & 0 & 1 & 1 \\
 1 & 1 & 1 & 0 & 0 & 0 \\
 0 & 0 & 0 & 1 & 1 & 1 \\
\end{psmallmatrix}
$ & $\begin{psmallmatrix}
 1 \\
 1 \\
 4 \\
 4 \\
 1 \\
 1 \\
\end{psmallmatrix}
$ & $\begin{psmallmatrix}
 0 & 0 & 1 & 1 & 2 & 2 \\
 0 & 0 & 1 & 1 & 2 & 2 \\
 3 & 3 & 0 & 4 & 4 & 4 \\
 4 & 4 & 4 & 0 & 3 & 3 \\
 2 & 2 & 1 & 1 & 0 & 0 \\
 2 & 2 & 1 & 1 & 0 & 0 \\
\end{psmallmatrix}
$\\\arrayrulecolor{white}
 \midrule
$\begin{psmallmatrix}
 -1 & 0 & 0 & 0 & 0 & 0 \\
 0 & -1 & 0 & 0 & 0 & 0 \\
 0 & 0 & -1 & 0 & 0 & 0 \\
 0 & 0 & 0 & -1 & 0 & 0 \\
 0 & 0 & 0 & 0 & -1 & 0 \\
 0 & 0 & 0 & 0 & 0 & -1 \\
 1 & 1 & 0 & 0 & -1 & -1 \\
 0 & 0 & 0 & 0 & 1 & -1 \\
 0 & 0 & 0 & 0 & -1 & 1 \\
 0 & 0 & 0 & 0 & 1 & 0 \\
 0 & 0 & 0 & 0 & 0 & 1 \\
 0 & 0 & 1 & 1 & 1 & 1 \\
\end{psmallmatrix}
$ & $\begin{psmallmatrix}
 1 \\
 1 \\
 5 \\
 5 \\
 1 \\
 1 \\
\end{psmallmatrix}
$ & $\begin{psmallmatrix}
 0 & 0 & 1 & 1 & 2 & 2 \\
 0 & 0 & 1 & 1 & 2 & 2 \\
 5 & 5 & 0 & 4 & 4 & 4 \\
 5 & 5 & 4 & 0 & 4 & 4 \\
 1 & 1 & 1 & 1 & 0 & 2 \\
 1 & 1 & 1 & 1 & 2 & 0 \\
\end{psmallmatrix}
$\\\arrayrulecolor{white}
 \midrule
$\begin{psmallmatrix}
 -1 & 0 & 0 & 0 & 0 & 0 \\
 0 & -1 & 0 & 0 & 0 & 0 \\
 0 & 0 & -1 & 0 & 0 & 0 \\
 0 & 0 & 0 & -1 & 0 & 0 \\
 0 & 0 & 0 & 0 & -1 & 0 \\
 0 & 0 & 0 & 0 & 0 & -1 \\
 1 & 1 & 0 & 0 & -1 & -1 \\
 0 & 0 & 0 & 0 & 1 & 0 \\
 0 & 0 & 0 & 0 & 0 & 1 \\
 0 & 0 & 1 & 1 & 1 & 1 \\
\end{psmallmatrix}
$ & $\begin{psmallmatrix}
 1 \\
 1 \\
 5 \\
 5 \\
 2 \\
 2 \\
\end{psmallmatrix}
$ & $\begin{psmallmatrix}
 0 & 0 & 1 & 1 & 2 & 2 \\
 0 & 0 & 1 & 1 & 2 & 2 \\
 5 & 5 & 0 & 4 & 4 & 4 \\
 5 & 5 & 4 & 0 & 4 & 4 \\
 2 & 2 & 2 & 2 & 0 & 2 \\
 2 & 2 & 2 & 2 & 2 & 0 \\
\end{psmallmatrix}
$\\\arrayrulecolor{white}
 \midrule
$\begin{psmallmatrix}
 -1 & 0 & 0 & 0 & 0 & 0 \\
 0 & -1 & 0 & 0 & 0 & 0 \\
 0 & 0 & -1 & 0 & 0 & 0 \\
 0 & 0 & 0 & -1 & 0 & 0 \\
 0 & 0 & 0 & 0 & -1 & 0 \\
 0 & 0 & 0 & 0 & 0 & -1 \\
 1 & 0 & 0 & 0 & 0 & -1 \\
 0 & 1 & 0 & 0 & 0 & -1 \\
 0 & 0 & 1 & 0 & 0 & -1 \\
 0 & 0 & -1 & 0 & 0 & 1 \\
 0 & 0 & 1 & 0 & 0 & 0 \\
 0 & 0 & 0 & 0 & 0 & 1 \\
 0 & 0 & 0 & 1 & 0 & 1 \\
 0 & 0 & 0 & 0 & 1 & 1 \\
\end{psmallmatrix}
$ & $\begin{psmallmatrix}
 1 \\
 1 \\
 1 \\
 3 \\
 3 \\
 1 \\
\end{psmallmatrix}
$ & $\begin{psmallmatrix}
 0 & 1 & 1 & 1 & 1 & 2 \\
 1 & 0 & 1 & 1 & 1 & 2 \\
 1 & 1 & 0 & 1 & 1 & 2 \\
 3 & 3 & 3 & 0 & 3 & 2 \\
 3 & 3 & 3 & 3 & 0 & 2 \\
 1 & 1 & 2 & 1 & 1 & 0 \\
\end{psmallmatrix}
$\\\arrayrulecolor{white}
 \midrule
$\begin{psmallmatrix}
 -1 & 0 & 0 & 0 & 0 & 0 \\
 0 & -1 & 0 & 0 & 0 & 0 \\
 0 & 0 & -1 & 0 & 0 & 0 \\
 0 & 0 & 0 & -1 & 0 & 0 \\
 0 & 0 & 0 & 0 & -1 & 0 \\
 0 & 0 & 0 & 0 & 0 & -1 \\
 1 & 0 & 0 & 0 & 0 & -1 \\
 0 & 1 & 0 & 0 & 0 & -1 \\
 0 & -1 & 1 & 0 & 0 & 0 \\
 0 & -1 & 0 & 0 & 0 & 1 \\
 0 & 1 & 0 & 0 & 0 & 0 \\
 0 & 0 & 0 & 0 & 0 & 1 \\
 0 & 1 & 0 & 1 & 0 & 0 \\
 0 & 0 & 0 & 0 & 1 & 1 \\
\end{psmallmatrix}
$ & $\begin{psmallmatrix}
 1 \\
 1 \\
 1 \\
 3 \\
 3 \\
 1 \\
\end{psmallmatrix}
$ & $\begin{psmallmatrix}
 0 & 1 & 1 & 1 & 1 & 2 \\
 1 & 0 & 1 & 1 & 1 & 2 \\
 1 & 2 & 0 & 1 & 1 & 1 \\
 3 & 2 & 3 & 0 & 3 & 3 \\
 3 & 3 & 3 & 3 & 0 & 2 \\
 1 & 2 & 1 & 1 & 1 & 0 \\
\end{psmallmatrix}$\\\arrayrulecolor{white}
 \midrule
$\begin{psmallmatrix}
 -1 & 0 & 0 & 0 & 0 & 0 \\
 0 & -1 & 0 & 0 & 0 & 0 \\
 0 & 0 & -1 & 0 & 0 & 0 \\
 0 & 0 & 0 & -1 & 0 & 0 \\
 0 & 0 & 0 & 0 & -1 & 0 \\
 0 & 0 & 0 & 0 & 0 & -1 \\
 1 & 0 & 0 & 0 & 0 & -1 \\
 -1 & 0 & 0 & 0 & 0 & 1 \\
 1 & 0 & 0 & 0 & 0 & 0 \\
 0 & 1 & 1 & 0 & 0 & -1 \\
 0 & 0 & 0 & 0 & 0 & 1 \\
 0 & 0 & 0 & 1 & 1 & 1 \\
\end{psmallmatrix}
$ & $\begin{psmallmatrix}
 1 \\
 2 \\
 2 \\
 4 \\
 4 \\
 1 \\
\end{psmallmatrix}$ & $\begin{psmallmatrix}
 0 & 1 & 1 & 1 & 1 & 2 \\
 2 & 0 & 1 & 2 & 2 & 3 \\
 2 & 1 & 0 & 2 & 2 & 3 \\
 4 & 4 & 4 & 0 & 3 & 3 \\
 4 & 4 & 4 & 3 & 0 & 3 \\
 2 & 1 & 1 & 1 & 1 & 0 \\
\end{psmallmatrix}$
\end{tabular}}
\caption{Undecomposable 6-dimensional smooth reflexive polytopes with barycenter at the origin (6-symplex  excluded).}
\label{table6}
\end{table}

\

\newpage 
\small{}

\end{document}